\documentclass[11pt,reqno]{article}
\usepackage{amsfonts}

\begin{document}

\def\R{{\mathbb R}}
\def\T{{\mathbb T}}
\def\S{{\mathbb S}}
\def\C{{\mathbb C}}
\def\Z{{\mathbb Z}}
\def\N{{\mathbb N}}
\def\H{{\mathbb H}}
\def\B{{\mathbb B}}
\def\diam{\mbox{\rm diam}}
\def\sn{\S^{n-1}}
\def\rr{{\cal R}}
\def\mt{{\Lambda}}
\def\e{\emptyset}
\def\dQ{\partial Q}
\def\dk{\partial K}
\def\endofproof{{\rule{6pt}{6pt}}}
\def\di{\displaystyle}
\def\dist{\mbox{\rm dist}}
\def\sa+{\Sigma_A^+}
\def\du{\frac{\partial}{\partial u}}
\def\dv{\frac{\partial}{\partial v}}
\def\dt{\frac{d}{d t}}
\def\dx{\frac{\partial}{\partial x}}
\def\con{\mbox{\rm const }}
\def\nn{{\cal N}}
\def\mm{{\cal M}}
\def\kk{{\cal K}}
\def\ll{{\cal L}}
\def\vv{{\cal V}}
\def\bb{{\cal B}}
\def\ma{\mm_{a}}
\def\lab{L_{ab}}
\def\mabn{\mm_{a}^N}
\def\man{\mm_a^N}
\def\labn{L_{ab}^N}
\def\fa{f^{(a)}}
\def\ff{{\cal F}}
\def\i{{\bf i}}
\def\gge{{\cal G}_\epsilon}
\def\gej{\chi^{(j)}_\mu}
\def\ge{\chi_\epsilon}
\def\geo{\chi^{(1)}_\mu}
\def\get{\chi^{(2)}_\mu}
\def\gei{\chi^{(i)}_{\mu}}
\def\gee{\chi_{\mu}}
\def\gett{\chi^{(2)}_{\mu}}
\def\geol{\chi^{(1)}_{\ell}}
\def\getl{\chi^{(2)}_{\ell}}
\def\geil{\chi^{(i)}_{\ell}}
\def\gee{\chi_{\ell}}
\def\tt{{\cal T}}
\def\uu{{\cal U}}
\def\wloc{W_{\epsilon}}
\def\Int{\mbox{\rm Int}}
\def\dist{\mbox{\rm dist}}
\def\pr{\mbox{\rm pr}}
\def\pp{{\cal P}}
\def\aa{{\cal A}}
\def\cc{{\cal C}}
\def\supp{\mbox{\rm supp}}
\def\Arg{\mbox{\rm Arg}}
\def\In{\mbox{\rm Int}}
\def\con{\mbox{\rm const}\;}
\def\Re{\mbox{\rm Re}}
\def\li{\mbox{\rm li}} 
\def\Seo{S^*_\epsilon(\Omega)}
\def\sdk{S^*_{\dk}(\Omega)}
\def\lae{\Lambda_{\epsilon}}
\def\ep{\epsilon}
\def\oo{{\cal O}}
\def\be{\begin{equation}}
\def\ee{\end{equation}}
\def\beqn{\begin{eqnarray*}}
\def\eeqn{\end{eqnarray*}}
\def\Pr{\mbox{\rm Pr}}

\def\gi{\gamma^{(i)}}
\def\ii{{\imath }}
\def\jj{{\jmath }}
\def\II{{\cal I}}
\def\ccij{ \cc_{i'_0,j'_0}[\eta]}
\def\dd{{\cal D}}
\def\la{\langle}
\def\ra{\rangle}
\def\bs{\bigskip}
\def\xio{\xi^{(0)}}
\def\xo{x^{(0)}}
\def\zo{z^{(0)}}
\def\Con{\mbox{\rm Const}\;}
\def\do{\partial \Omega}
\def\dk{\partial K}
\def\dl{\partial L}
\def\ll{{\cal L}}
\def\kk{{\cal K}}
\def\kk{{\cal K}}
\def\pr{{\rm pr}}
\def\ff{{\cal F}}
\def\G{{\cal G}}
\def\C{{\bf C}}
\def\dist{{\rm dist}}
\def\dds{\frac{d}{ds}}
\def\con{{\rm const}\;}
\def\Con{{\rm Const}\;}
\def\di{\displaystyle}
\def\oo{\mbox{\rm O}}
\def\hess{\mbox{\rm Hess}}
\def\gi{\gamma^{(i)}}
\def\endofproof{{\rule{6pt}{6pt}}}
\def\xm{x^{(m)}}
\def\vm{\varphi^{(m)}}
\def\km{k^{(m)}}
\def\dm{d^{(m)}}
\def\kam{\kappa^{(m)}}
\def\dem{\delta^{(m)}}
\def\xim{\xi^{(m)}}
\def\ep{\epsilon}
\def\ms{\medskip}
\def\ex{\mbox{\rm extd}}

% some new definitions
\def\clip{C^{\mbox{\footnotesize \rm Lip}}}
\def\wlocs{W^s_{\mbox{\footnote\rm loc}}}
\def\Lip{\mbox{\rm Lip}}

\def\Xr{X^{(r)}}
\def\lip{\mbox{{\footnotesize\rm Lip}}}
\def\Vol{\mbox{\rm Vol}}

\def\naf{\nabla f(z)}
\def\so{\sigma_0}
\def\Xo{X^{(0)}}
\def\z1{z^{(1)}}
\def\Vo{V^{(0)}}
\def\Yo{Y{(0)}}

\def\uo{u^{(0)}}
\def\vo{v^{(0)}}
\def\no{\nu^{(0)}}
\def\psa{\partial^{(s)}_a}
\def\hcd{\hc^{(\delta)}}
\def\Md{M^{(\delta)}}
\def\Uo{U^{(1)}}
\def\Ut{U^{(2)}}
\def\Uj{U^{(j)}}
\def\no{n^{(1)}}
\def\nt{n^{(2)}}
\def\nj{n^{(j)}}
\def\ccm{\cc^{(m)}}

\def\Ul{U^{(\ell)}}
\def\Uj{U^{(j)}}
\def\wl{w^{(\ell)}}
\def\Vl{V^{(\ell)}}
\def\Ujj{U^{(j+1)}}
\def\wjj{w^{(j+1)}}
\def\Vjj{V^{(j+1)}}
\def\Ujo{U^{(j_0)}}
\def\wjo{w^{(j_0)}}
\def\Vjo{V^{(j_0)}}
\def\vj{v^{(j)}}
\def\vl{v^{(\ell)}}

\def\f0{f^{(0)}}

\def\gl{\gamma_\ell}
\def\id{\mbox{\rm id}}
\def\piU{\pi^{(U)}}

\def\cca{C^{(a)}}
\def\bba{B^{(a)}}
\def\co{\; \stackrel{\circ}{C}}
\def\Lo{\; \stackrel{\circ}{L}}

\def\oV{\overline{V}}
\def\saa{\Sigma^+_A}
\def\sa{\Sigma_A}
\def\mta{\Lambda(A, \tau)}
\def\mtaa{\Lambda^+(A, \tau)}

\def\Int{\mbox{\rm Int}}
\def\epo{\ep^{(0)}}
\def\pH{\partial \H^{n+1}}
\def\sh{S^*(\H^{n+1})}
\def\zoo{z^{(1)}}
\def\yoo{y^{(1)}}
\def\xoo{x^{(1)}}

%balls2.tex

%General
\def\supp{\mbox{\rm supp}}
\def\Arg{\mbox{\rm Arg}}
\def\In{\mbox{\rm Int}}
\def\diam{\mbox{\rm diam}}
\def\e{\emptyset}
\def\endofproof{{\rule{6pt}{6pt}}}
\def\di{\displaystyle}
\def\dist{\mbox{\rm dist}}
\def\con{\mbox{\rm const }}
\def\Box{\spadesuit}
\def\Int{\mbox{\rm Int}}
\def\dist{\mbox{\rm dist}}
\def\pr{\mbox{\rm pr}}
\def\be{\begin{equation}}
\def\ee{\end{equation}}
\def\beqn{\begin{eqnarray*}}
\def\eeqn{\end{eqnarray*}}
\def\la{\langle}
\def\ra{\rangle}
\def\bs{\bigskip}
\def\Con{\mbox{\rm Const}\;}
\def\clip{C^{\mbox{\footnotesize \rm Lip}}}
\def\wlocs{W^s_{\mbox{\footnote\rm loc}}}
\def\Lip{\mbox{\rm Lip}}
\def\lip{\mbox{\footnotesize\rm Lip}}
\def\Re{\mbox{\rm Re}}
\def\li{\mbox{\rm li}} 
\def\ep{\epsilon}
\def\ms{\medskip}
\def\dds{\frac{d}{ds}}
\def\oo{\mbox{\rm O}}
\def\hess{\mbox{\rm Hess}}
\def\id{\mbox{\rm id}}
\def\ii{{\imath }}
\def\jj{{\jmath }}
\def\graph{\mbox{\rm graph}}
\def\span{\mbox{\rm span}}

% Boldface
\def\i{{\bf i}}
\def\C{{\bf C}}

% Caligraphic
\def\ss{{\cal S}}
\def\tt{{\cal T}}
\def\E{{\cal E}}
\def\B{{\cal B}}
\def\rr{{\cal R}}
\def\nn{{\cal N}}
\def\mm{{\cal M}}
\def\kk{{\cal K}}
\def\ll{{\cal L}}
\def\vv{{\cal V}}
\def\ff{{\cal F}}
\def\hh{{\cal H}}
\def\tt{{\cal T}}
\def\uu{{\cal U}}
\def\cc{{\cal C}}
\def\pp{{\cal P}}
\def\aa{{\cal A}}
\def\oo{{\cal O}}
\def\II{{\cal I}}
\def\dd{{\cal D}}
\def\ll{{\cal L}}
\def\ff{{\cal F}}
\def\G{{\cal G}}
\def\W{{\cal W}}
\def\hhs{\hh^s}
\def\thhs{\widetilde{\hh}^s}
\def\hhhs{\widehat{\hh}^s}

% Hat
\def\hs{\hat{s}}
\def\hz{\hat{z}}
\def\hL{\hat{L}}
\def\hl{\hat{l}}
\def\hl{\hat{l}}
\def\hc{\hat{\cc}}
\def\hbb{\widehat{\cal B}}
\def\hu{\hat{u}}
\def\hX{\hat{X}}
\def\hx{\hat{x}}
\def\hu{\hat{u}}
\def\hv{\hat{v}}
\def\hQ{\hat{Q}}
\def\hC{\widehat{C}}
\def\hF{\hat{F}}
\def\hf{\hat{f}}
\def\hii{\hat{\ii}}
\def\hr{\hat{r}}
\def\hq{\hat{q}}
\def\hy{\hat{y}}
\def\hZ{\widehat{Z}}
\def\hz{\hat{z}}
\def\hE{\widehat{E}}
\def\hR{\widehat{R}}
\def\hell{\hat{\ell}}
\def\hs{\hat{s}}
\def\hW{\widehat{W}}
\def\hS{\widehat{S}}
\def\hV{\widehat{V}}
\def\hB{\widehat{B}}
\def\hhh{\widehat{\cal H}}
\def\hK{\widehat{K}}
\def\hU{\widehat{U}}
\def\hhh{\widehat{\hh}}
\def\hdd{\widehat{\dd}}
\def\hZ{\widehat{Z}}

\def\hal{\hat{\alpha}}
\def\hbe{\hat{\beta}}
\def\hg{\hat{\gamma}}
\def\hrho{\hat{\rho}}
\def\hd{\hat{\delta}}
\def\hphi{\hat{\phi}}
\def\hmu{\hat{\mu}}
\def\hnu{\hat{\nu}}
\def\hsi{\hat{\sigma}}
\def\htau{\hat{\tau}}
\def\hpi{\hat{\pi}}
\def\hep{\hat{\epsilon}}
\def\hxi{\hat{\xi}}
\def\hLa{\widehat{\Lambda}^u}
\def\hPhi{\widehat{\Phi}}
\def\hPsi{\widehat{\Psi}}
\def\hPhii{\widehat{\Phi}^{(i)}}

% Tilde
\def\tc{\tilde{C}}
\def\tg{\tilde{\gamma}}  
\def\tV{\widetilde{V}}
\def\tC{\widetilde{\cc}}
\def\tr{\tilde{R}}
\def\tb{\tilde{b}}
\def\tt{\tilde{t}}
\def\tx{\tilde{x}}
\def\tp{\tilde{p}}
\def\tz{\tilde{Z}}
\def\tZ{\tilde{Z}}
\def\tF{\tilde{F}}
\def\tf{\tilde{f}}
\def\tp{\tilde{p}}
\def\te{\tilde{e}}
\def\tv{\tilde{v}}
\def\tu{\tilde{u}}
\def\tw{\tilde{w}}
\def\ts{\tilde{\sigma}}
\def\tr{\tilde{r}}
\def\tS{\tilde{S}}
\def\tP{\widetilde{\Pi}}
\def\ttau{\tilde{\tau}}
\def\tLip{\widetilde{\Lip}}
\def\tz{\tilde{z}}
\def\tS{\tilde{S}}
\def\tts{\tilde{\sigma}}
\def\tVl{\widetilde{V}^{(\ell)}}
\def\tVj{\widetilde{V}^{(j)}}
\def\tVo{\widetilde{V}^{(1)}}
\def\tVj{\widetilde{V}^{(j)}}
\def\tPsi{\tilde{\Psi}}
 \def\tp{\tilde{p}}
 \def\tVjo{\widetilde{V}^{(j_0)}}
\def\tvj{\tilde{v}^{(j)}}
\def\tVjj{\widetilde{V}^{(j+1)}}
\def\tvl{\tilde{v}^{(\ell)}}
\def\tVt{\widetilde{V}^{(2)}}
\def\tR{\tilde{R}}
\def\tQ{\tilde{Q}}
\def\oL{\tilde{\Lambda}}
\def\tq{\tilde{q}}
\def\tx{\tilde{x}}
\def\ty{\tilde{y}}
\def\tz{\tilde{z}}
\def\txo{\tilde{x}^{(0)}}
\def\tso{\tilde{\sigma}_0}
\def\tmt{\tilde{\Lambda}}
\def\tg{\tilde{g}}
\def\tsi{\tilde{\sigma}}
\def\ttt{\tilde{t}}
\def\tC{\tilde{C}}
\def\tc{\tilde{c}}
\def\tell{\tilde{\ell}}
\def\trho{\tilde{\rho}}
\def\ts{\tilde{s}}
\def\tB{\widetilde{B}}
\def\thh{\widetilde{\cal H}}
\def\tV{\widetilde{V}}
\def\trr{\tilde{r}}
\def\te{\tilde{e}}
\def\tv{\tilde{v}}
\def\tu{\tilde{u}}
\def\tw{\tilde{w}}
\def\trho{\tilde{\rho}}
\def\tell{\tilde{\ell}}
\def\tz{\tilde{Z}}
\def\tF{\tilde{F}}
\def\tf{\tilde{f}}
\def\tp{\tilde{p}}
\def\ttau{\tilde{\tau}}
\def\tz{\tilde{z}}
\def\tg{\tilde{\gamma}}  
\def\tV{\widetilde{V}}
\def\tC{\widetilde{\cc}}
\def\tLa{\widetilde{\Lambda}^u}
\def\tR{\widetilde{R}}
\def\tr{\tilde{r}}
\def\tc{\widetilde{C}}
\def\tD{\widetilde{D}}
\def\tt{\tilde{t}}
\def\tx{\tilde{x}}
\def\tp{\tilde{p}}
\def\tS{\tilde{S}}
\def\tts{\tilde{\sigma}}
\def\tZ{\widetilde{Z}}
\def\tdelta{\tilde{\delta}}
\def\th{\tilde{h}}
\def\tB{\widetilde{B}}
\def\thh{\widetilde{\hh}}
\def\tep{\tilde{\ep}}
\def\tE{\widetilde{E}}
\def\tu{\tilde{u}}
\def\txi{\tilde{\xi}}
\def\teta{\tilde{\eta}}

% Specific
\def\sr{{\sc r}}
\def\mt{{\Lambda}}
\def\do{\partial \Omega}
\def\dk{\partial K}
\def\dl{\partial L}
\def\wloc{W_{\epsilon}}
\def\piU{\pi^{(U)}}
\def\Rio{\R_{i_0}}
\def\Ri{\R_{i}}
\def\Rii{\R^{(i)}}
\def\Riii{\R^{(i-1)}}
\def\hRii{\widehat{\R}_i}
\def\hRiio{\widehat{\R}_{(i_0)}}
\def\Eii{E^{(i)}}
\def\Eio{E^{(i_0)}}
\def\Rj{\R_{j}}
\def\Vio{{\cal V}^{i_0}}
\def\Vi{{\cal V}^{i}}
\def\Wio{W^{i_0}}
\def\Wioo{W^{i_0-1}}
\def\hi{h^{(i)}}
\def\Psii{\Psi^{(i)}}
\def\pii{\pi^{(i)}}
\def\piii{\pi^{(i-1)}}
\def\gxyii{g_{x,y}^{i-1}}
\def\span{\mbox{\rm span}}
\def\Jac{\mbox{\rm Jac}}
\def\Vol{\mbox{\rm Vol}}
\def\limp{\lim_{p\to\infty}}
\def\hh{{\mathcal H}}

\def\yijl{Y_{i,j}^{(\ell)}}
\def\xijl{X_{i,j}^{(\ell)}}
\def\hyijl{\widehat{Y}_{i,j}^{(\ell)}}
\def\hxijl{\widehat{X}_{i,j}^{(\ell)}}
\def\eijl{\eta_{i,j}^{(\ell)}}
\def\J{\sf J}
\def\Gl{\Gamma_\ell}

\def\hLao{\widehat{\Lambda}^{u,1}}
\def\tLao{\widetilde{\Lambda}^{u,1}}
\def\Lao{\Lambda^{u,1}}
\def\cLao{\check{\Lambda}^{u,1}}
\def\cB{\check{B}}
\def\tpi{\tilde{\pi}}
\def\hcc{\widehat{\cc}}

\def\tE{\widetilde{E}}
\def\Wuo{W^{u,1}}
\def\Wut{W^{u,2}}

\def\wo{w^{(1)}}
\def\vo{v^{(1)}}
\def\uo{u^{(1)}}
\def\wt{w^{(2)}}
\def\vt{v^{(2)}}
\def\ut{u^{(2)}}
\def\Wo{W^{(1)}}
\def\Vo{V^{(1)}}
\def\Uo{U^{(1)}}
\def\Wt{W^{(2)}}
\def\Vt{V^{(2)}}
\def\Ut{U^{(2)}}
\def\tmu{\tilde{\mu}}
\def\tla{\tilde{\lambda}}
\def\diamf{{\rm\footnotesize diam}}
\def\Intu{\mbox{\rm Int}^u}
\def\Ints{\mbox{\rm Int}^s}

\begin{center}
%\noindent
{\Large\bf  Regular decay of ball diameters\\
and spectra of Ruelle operators for contact Anosov flows}
\end{center}

\begin{center}
{\sc Luchezar Stoyanov}
\end{center}

\bs

%\footnotesize

\noindent
{\it Abstract.} For Anosov flows on compact Riemann manifolds we study the rate of
decay along the flow of diameters of balls $B^s(x,\ep)$ on local stable manifolds at 
Lyapunov regular points $x$. We prove that this decay rate is similar for all
sufficiently small values of $\ep > 0$. From this and the main result in \cite{kn:St1},
we derive strong spectral estimates for Ruelle transfer operators
for contact Anosov flows with Lipschitz local stable holonomy maps. These apply
in particular to geodesic flows on compact locally symmetric manifolds of strictly
negative curvature.
As is now well known, such spectral estimates have deep implications in some related 
areas, e.g. in studying analytic properties of  Ruelle zeta functions and partial differential 
operators,  asymptotics of closed orbit counting functions, etc.

\section{\sc Introduction}
\renewcommand{\theequation}{\arabic{section}.\arabic{equation}}

Consider a non-linear system of differential equations of the form
\be
\dot{x}(t) = f(t,x)\;,
\ee
where $f : U\times \R \longrightarrow \R^n$ is a continuously differentiable map for some 
open ball $U$ with center $0$ in $\R^n$ and $f(t,0) = 0$ for all $t \in \R$. Assuming that the null
solution of (1.1) is asymptotically stable (see e.g. \cite{kn:CL}), one defines a semi-flow
$\varphi_t : U\times [0,\infty) \longrightarrow U$ such that $\varphi_t(z,s) = x(t+s)$, where
$x$ is the solution of (1.1) with $x(s) = z$. One may then ask the question whether for all
sufficiently small $0 < \delta_1 < \delta_2$ there exists a constant $C > 0$ depending only
on $\delta_1$ and $\delta_2$ (and $f$) such that 
$\diam(\varphi_t(B(0,\delta_2))) \leq C \diam(\varphi_t(B(0,\delta_1)))$ for all $t \geq 0$,
where $B(0,\delta)$ denotes the (closed) ball with center $0$ and radius $\delta$ in $\R^n$.
We do not know what happens in the general case, however it follows from the 
arguments in the present paper that under a certain (Lyapunov regularity) condition at $0$, the 
answer to the above question is affirmative.

In fact,  we consider a more complicated situation.
Let $\phi_t : M \longrightarrow M$ be a $C^2$ Anosov flow on a $C^2$ compact Riemann manifold $M$.
%In most of this paper we will assume that the local stable holonomy maps associated with the flow
%are Lipschitz.
For any $x \in M$ and  a sufficiently small $\delta > 0$ consider the closed $\delta$-ball
$B^s (x,\delta) = \{ y\in W^s_{\ep}(x) : d(x, y) \leq \delta  \}\;$
on the local stable manifold $W^s_\ep(x)$. For any $y\in W^s_\ep(x)$ we know that 
$d(\phi_t(x),\phi_t(y)) \to 0$ exponentially fast as $t \to \infty$. Moreover, we have uniform estimates
for the exponential rate of convergence, so for any $x \in M$ and any given $\delta > 0$, 
$\diam(\phi_t(B^s(x,\delta))) \to 0$ exponentially fast as $t \to \infty$. However, in general it is 
not clear whether for any constants $0 < \delta_1 < \delta_2$ the ratio 
$$\frac{\diamf(\phi_t(B^s(x,\delta_2)))}{\diamf(\phi_t(B^s(x,\delta_1)))}$$
is uniformly bounded  for $t > 0$ and $x\in M$ (although a similar property is obviously 
satisfied by the linearized flow $d\phi_t$, considering balls on corresponding tangent planes).
It appears that in general this problem is rather subtle, and  it is not clear at all whether one 
should expect a positive solution  without any extra assumptions.

Here we consider a similar problem on the set $\ll$ of Lyapunov regular points in $M$ --
see section 3.1 for the terminology. We prove the following.

\bs

\noindent
{\bf Theorem 1.1.} {\it  For every $\ep > 0$ there exist Lyapunov $\ep$-regularity functions
$\omega : \ll \longrightarrow (0,1)$ and $G: \ll \longrightarrow [1,\infty)$ such that
for  any $0 < \delta_1 < \delta_2$ there exists a constant $K = K(\delta_1,\delta_2) \geq 1$
with 
$$\diam(\phi_t(B^s(x, \delta_2)) \leq K G(x) \diam (\phi_t(B^s(x,\delta_1))$$
for all  $x\in \ll$ with $\delta_2 \leq \omega(x)$ and all $ t > 0$.}

\ms

A similar result can be proved for non-uniformly hyperbolic flows.

The above has an important consequence concerning cylinders in a symbolic coding of the flow 
defined by  means of a Markov family -- see Theorem 4.2 below for details.

The motivation for this work came from \cite{kn:St1}, where assuming the properties (a) and (b) in
Theorem 4.2 below, Lipschitzness of the local stable holonomy  maps and a certain 
non-integrability condition 
we prove strong spectral estimates for arbitrary potentials over basic sets for Axiom A flows, similar to 
those established by Dolgopyat \cite{kn:D} for geodesic flows on compact surfaces (for general potentials)
and transitive Anosov flows on compact manifolds  with $C^1$ jointly non-integrable 
horocycle foliations (for the Sinai-Bowen-Ruelle potential).
It is known that such strong spectral estimates  lead to deep results in a variety of 
areas which are difficult  (if not impossible) to obtain by  other means (see e.g. \cite{kn:PoS1},
\cite{kn:PoS2}, \cite{kn:PoS3}, \cite{kn:An},  \cite{kn:PeS1} \cite{kn:PeS2}, \cite{kn:PeS3}).

Let $\rr = \{R_i\}_{i=1}^k$ be a Markov family for $\phi_t$  consisting of 
rectangles $R_i = [U_i ,S_i ]$, where $U_i$ (resp. $S_i$) are (admissible) subsets of 
$W^u_{\ep}(z_i)$ (resp. $W^s_{\ep}(z_i)$) for some $\ep > 0$ and $z_i\in M$ 
(cf. section 2 for details). The {\it first return time function}  
$\tau : R = \cup_{i=1}^k R_i \longrightarrow [0,\infty)$ 
and the standard Poincar\'e map $\pp: R \longrightarrow R$ are then naturally defined.
Setting $U = \cup_{i=1}^k U_i$, the {\it shift map} $\sigma : U \longrightarrow U$
defined by $\sigma = \piU\circ \pp$, where $\piU : R \longrightarrow U$
is the projection along the leaves of local stable manifolds, provides
a natural symbolic coding of the flow. 
%and therefore defines a family of basis sets for the topology
%in the corresponding symbol space -- the so called cylinders. 
To avoid dealing with boundary points in $U$, consider the set $\hU$ of all 
$u \in U$ whose orbits do not have common points with the boundary of $R$ 
(see section 2). 

%Given an admissible finite string $\ii = (i_0,i_1, \ldots,i_m)$  of integers $i_j \in \{ 1, \ldots,k\}$
%(see section 2) denote by  $C[\ii]$ the closure (in $U$) of the set of those $x\in \hU$ so that 
%$\sigma^j(x) \in U_{i_j}$ for all $j = 0,1, \ldots,m$, is called a {\it cylinder} of length $m$  in $U$.

Given a Lipschitz real-valued function $f$  on $\hU$, set $g = g_f = f - P\tau$, where 
$P = P_f\in \R$ is the unique 
number such that the topological pressure $\Pr_\sigma(g)$ of $g$ with respect to $\sigma$ is 
zero (cf. e.g. \cite{kn:PP}). For $a, b\in \R$, one defines the {\it Ruelle transfer operator}
$L_{g-(a+\i b)\tau} : \clip (\hU) \longrightarrow \clip (\hU)$ in the usual way (cf. section 2). Here
$\clip (\hU)$ is the space of Lipschitz functions $g: \hU \longrightarrow \C$. By 
$\Lip(g)$ we denote the Lipschitz constant of $g$ and  by $\| g\|_0$ the {\it standard $\sup$ norm}  
of $g$ on $\hU$.

We will say  that the {\it Ruelle transfer operators related to the function $f$ on $U$ are 
eventually contracting}  if for every $\epsilon > 0$ there exist constants $0 < \rho < 1$, $a_0 > 0$ 
and  $C > 0$ such  that if $a,b\in \R$ satisfy $|a| \leq a_0$ and $|b| \geq 1/a_0$, then for every integer 
$m > 0$ and every  $h\in \clip (\hU)$ we have
$$\|L_{f -(P_f+a+ \i b)\tau}^m h \|_{\lip,b} \leq C \;\rho^m \;|b|^{\ep}\; \| h\|_{\lip,b}\; ,$$
where the norm $\|.\|_{\lip,b}$ on $\clip (\hU)$ is defined by 
$\| h\|_{\lip,b} = \|h\|_0 + \frac{\lip(h)}{|b|}$.
This implies in particular that the spectral radius  of $L_{f-(P_f+ a+\i b)\tau}$ on $\clip(\hU)$ 
does not exceed  $\rho$.

From Theorem 1.1 (or rather its consequence -- Theorem 4.2 below) and the main 
result in \cite{kn:St1}  we derive the following.

\bs

\noindent
{\bf Theorem 1.2.} {\it  Let $\phi_t : M \longrightarrow M$ be a $C^2$ transitive contact Anosov flow on
a $C^2$ compact  Riemann manifold with uniformly Lipschitz local stable holonomy maps.
Then for any  Lipschitz real-valued function 
$f$  on $U$ the Ruelle transfer operators related to  $f$  are eventually contracting}.

\bs

The reader is referred to section 2 below for the definition of local holonomy maps. 
In general these  are only H\"older continuous. It is known that uniform
Lipschitzness of the local stable holonomy maps can be derived from certain bunching 
condition concerning the rates of expansion/contraction of the flow along local unstable/stable 
manifolds over $M$ (see \cite{kn:Ha},  \cite{kn:PSW}).

A result similar to Theorem 1.2 is true for general (non necessarily contact) Anosov flows, however one has to
assume in addition a local non-integrability condition (see condition (LNIC) in \cite{kn:St1}).
Using a smoothing procedure as in \cite{kn:D}, an estimate similar to that in Theorem 1.2  
holds for the Ruelle operator acting on the space $\ff_\gamma(U)$ of H\"older continuous functions 
with respect to an appropriate norm.

For geodesic flows  on locally symmetric spaces of negative curvature it is well known that
the local stable and unstable manifolds are smooth ($C^{\infty}$), so the corresponding local holonomy 
maps are smooth as well. 
%(see e.g. Theorem 1 and fact (2) on p. 647 in \cite{kn:Ha}).
Thus, as an immediate consequence of Theorem 1.2 one obtains  the following.

\bs

\noindent
{\bf Theorem 1.3.} {\it Let $X$ be a compact locally symmetric space of negative curvature 
and let\\ $\phi_t : M = S^*(X) \longrightarrow M$ be the geodesic flow on $X$. Then for any  
Lipschitz real-valued function $f$  on $U$ the Ruelle transfer operators related to  $f$  are 
eventually contracting}.

\bs

As mentioned above, there are various consequences that can be derived from results like
Theorem 1.2 (or Theorem 1.3).  Here we state one of these.

As in \cite{kn:St1}, one can use Theorem 1.2  and an argument of Pollicott and Sharp \cite{kn:PoS1} 
to get certain information  about the {\it Ruelle zeta function}
$$\zeta(s) = \prod_{\gamma} (1- e^{-s\ell(\gamma)})^{-1}\;,$$
where $\gamma$ runs over the set of primitive  closed orbits of $\phi_t$
and $\ell(\gamma)$ is the least period of $\gamma$.  Let $h_T$ denote the
 {\it topological entropy} of $\phi_t$.\\

\noindent
{\bf Corollary 1.4.} {\it Under the assumptions in Theorems 1.2 or 1.3, 
the  zeta function $\zeta(s)$ of the flow $\phi_t: M \longrightarrow M$ has an analytic  
and non-vanishing continuation in a half-plane $\Re(s) > c_0$ for some $c_0 < h_T$ except for a 
simple pole at $s = h_T$.  Moreover, there exists $c \in (0, h_T)$ such that
$$\pi(\lambda) = \# \{ \gamma : \ell(\gamma) \leq \lambda\} = \li(e^{h_T \lambda}) + O(e^{c\lambda})$$
as $\lambda\to \infty$, where $\di \li(x) = \int_2^x \frac{du}{\log u} \sim \frac{x}{\log x}$ as 
$x \to \infty$. }

\bigskip

In fact, a direct application of Theorem 5 in \cite{kn:PeS3} gives a more precise estimate of the
number of closed trajectories of the flow with primitive periods lying in exponentially shrinking intervals 
-- we refer the reader to section 6 in \cite{kn:PeS3} for details.

Section 2 contains some basic definitions and preliminary facts. In section 3 we compare diameters of balls
with respect to Bowen's metric on unstable manifolds and prove the analogue of Theorem 1.1 for unstable 
manifolds. From this Theorem 1.1 is derived easily. Finally in section 4 we consider cylinders in the set $U$
defined by means of a Markov family, and prove two properties of the decay rates of the diameters of such
cylinders (Theorem 4.2), assuming that the local stable holonomy maps are uniformly Lispchitz.
We do not know whether the same properties hold for any Anosov flow. Theorem 1.2 is then derived
using Theorem 4.2 and the argument in section 5 of \cite{kn:St1}.

%\newpage

%Sect. 2
\section{\sc Preliminaries}
\setcounter{equation}{0}

%\subsection{\sc Basics}

Throughout this paper $M$ denotes a $C^2$ compact Riemann manifold,  and 
$\phi_t : M \longrightarrow M$ ($t\in \R$) a $C^2$ flow on $M$. The flow 
is called {\it hyperbolic} if $M$ contains no fixed points  and there exist  constants 
$C > 0$ and $0 < \lambda < 1$ such that there exists a $d\phi_t$-invariant decomposition 
$T_xM = E^0(x) \oplus E^u(x) \oplus E^s(x)$ of $T_xM$ ($x \in M$) into a 
direct sum of non-zero linear subspaces, where $E^0(x)$ is the one-dimensional subspace 
determined by the direction of the flow at $x$, $\| d\phi_t(u)\| \leq C\, \lambda^t\, \|u\|$ for all  
$u\in E^s(x)$ and $t\geq 0$, and $\| d\phi_t(u)\| \leq C\, \lambda^{-t}\, \|u\|$ for all $u\in E^u(x)$ 
and  $t\leq 0$. The flow $\phi_t$ is called an {\it Anosov flow} on $M$ if 
the periodic points are dense in $M$ (see e.g. \cite{kn:KH}). 
The flow is called {\it transitive} if it has a dense orbit, and {\it contact} if there exists a $C^2$ flow 
invariant one form $\omega$ on $M$ such that $\omega \wedge (d\omega)^n$ is nowhere zero,
where $\dim(M) = 2n+1$.

For $x\in M$ and a sufficiently small $\epsilon > 0$ let 
$$\wloc^s(x) = \{ y\in M : d (\phi_t(x),\phi_t(y)) \leq \epsilon \: \mbox{\rm for all }
\: t \geq 0 \; , \: d (\phi_t(x),\phi_t(y)) \to_{t\to \infty} 0\: \}\; ,$$
$$\wloc^u(x) = \{ y\in M : d (\phi_t(x),\phi_t(y)) \leq \epsilon \: \mbox{\rm for all }
\: t \leq 0 \; , \: d (\phi_t(x),\phi_t(y)) \to_{t\to -\infty} 0\: \}$$
be the (strong) {\it stable} and {\it unstable manifolds} of size $\epsilon$. Then
$E^u(x) = T_x \wloc^u(x)$ and $E^s(x) = T_x \wloc^s(x)$. 
Given $\delta > 0$, set $E^u(x;\delta) = \{ u\in E^u(x) : \|u\| \leq \delta\}$;
$E^s(x;\delta)$ is defined similarly. 
%For any $A\subset M$ and  $I \subset \R$  denote
%$\phi_I(A) = \{\; \phi_t(y)\; : \; y\in A, t \in I \; \}.$

It follows from the hyperbolicity of the flow  that if  $\epsilon_0 > 0$ is sufficiently small,
there exists $\ep_1 > 0$ such that if $x,y\in M$ and $d (x,y) < \ep_1$, 
then $W^s_{\ep_0}(x)$ and $\phi_{[-\ep_0,\ep_0]}(W^u_{\ep_0}(y))$ intersect at exactly 
one point $[x,y ] \in M$  (cf. \cite{kn:KH}). That is, there exists a unique 
$t\in [-\ep_0, \ep_0]$ such that $\phi_t([x,y]) \in W^u_{\ep_0}(y)$. 
%Setting $\Delta(x,y) = t$, defines the so called {\it temporal distance function}. 
For $x, y\in M$ with $d (x,y) < \ep_1$, define
$\pi_y(x) = [x,y] = W^s_{\ep}(x) \cap \phi_{[-\ep_0,\ep_0]} (W^u_{\ep_0}(y))\;.$
Thus, for a fixed $y \in M$, $\pi_y : W \longrightarrow \phi_{[-\ep_0,\ep_0]} (W^u_{\ep_0}(y))$ is the
{\it projection} along local stable manifolds defined on a small open neighbourhood $W$ of $y$ in $M$.
Choosing $\ep_1 \in (0,\ep_0)$ sufficiently small,  the restriction 
$\pi_y: \phi_{[-\ep_1,\ep_1]} (W^u_{\ep_1}(x)) \longrightarrow \phi_{[-\ep_0,\ep_0]} (W^u_{\ep_0}(y))$
is called a local stable holonomy map\footnote{In a similar way one can define
holonomy maps between any two sufficiently close local transversals to stable laminations; see e.g.
\cite{kn:PSW}.}. 
%Combining it with a shift along the flow we get another {\it local stable holonomy  map}
%$\hhs_{x,y} : W^u_{\ep_1}(x) \cap \mt  \longrightarrow W^u_{\ep_0}(y) \cap \mt$.
%In a similar way one defines local holonomy maps along unstable laminations.

We will say that $A$ is an {\it admissible subset} of $W^u_{\ep}(z)$ ($z\in M$) 
if $A$ coincides with the closure of its interior in $W^u_\ep(z)$. Admissible subsets of
$W^s_\ep(z)$ are defined similarly. As in \cite{kn:D}, a subset $R$ of $\mt$ will be called a {\it rectangle} 
if it has the form $R = [U,S] = \{ [x,y] : x\in U, y\in S\}$, where $U$ and $S$ are admissible 
subsets of $W^u_\ep(z)$ and $W^s_\ep(z)$, respectively, for some $z\in M$. 
In what follows we will denote by  $\Intu(U)$ the {\it interior} of $U$ in the set $\wloc^u(z)$. 
In a similar way we  define $\Ints(S)$, and then set $\Int(R) = [\Intu(U), \Ints(S)]$.
Given $\xi = [x,y] \in R$, set $W^u_R(\xi) = [U,y] = \{ [x',y] : x'\in U\}$ and
$W^s_R(\xi) = [x,S] = \{[x,y'] : y'\in S\} \subset W^s_{\ep_0}(x)$. The {\it interiors}
of these sets in the corresponding leaves are defined by $\Intu(W^u_R(\xi)) = [\Intu(U),y]$
and $\Ints(W^s_R(\xi)) = [x,\Ints(S)]$.

Let $\rr = \{ R_i\}_{i=1}^k$ be a family of rectangles with $R_i = [U_i  , S_i ]$,
$U_i \subset \wloc^u(z_i)$ and $S_i \subset \wloc^s(z_i)$, respectively, 
for some $z_i\in M$.  Since the set $\ll$ of Lyapunov regular points (see section 3.1
below) is dense, without loss of generality we will assume that $z_i \in R_i \cap \ll$.
Set $R =  \cup_{i=1}^k R_i$.
The family $\rr$ is called {\it complete} if  there exists $T > 0$ such that for every $x \in M$,
$\phi_{t}(x) \in R$ for some  $t \in (0,T]$.  Given a complete family $\rr$, the related
{\it Poincar\'e map} $\pp: R \longrightarrow R$
is defined by $\pp(x) = \phi_{\tau(x)}(x) \in R$, where
$\tau(x) > 0$ is the smallest positive time with $\phi_{\tau(x)}(x) \in R$.
The function $\tau$  is called the {\it first return time}  associated with $\rr$. 
A complete family $\rr = \{ R_i\}_{i=1}^k$ of rectangles in $M$ is called a 
{\it Markov family} of size $\chi > 0$ for the  flow $\phi_t$ if $\diam(R_i) < \chi$ for all $i$ and: 
(a)  for any $i\neq j$ and any $x\in \Int(R_i) \cap \pp^{-1}(\Int(R_j))$ we have   
$\pp(\Ints (W_{R_i}^s(x)) ) \subset \Ints (W_{R_j}^s(\pp(x)))$ and 
$\pp(\Intu(W_{R_i}^u(x))) \supset \Intu(W_{R_j}^u(\pp(x)))$;
(b) for any $i\neq j$ at least one of the sets $R_i \cap \phi_{[0,\chi]}(R_j)$ and
$R_j \cap \phi_{[0,\chi]}(R_i)$ is empty.

The existence of a Markov family $\rr$ of an arbitrarily small size $\chi > 0$ for $\phi_t$
follows from the construction of Bowen \cite{kn:B} (cf. also  Ratner \cite{kn:Ra}). 

From now on we will assume that $\rr = \{ R_i\}_{i=1}^k$ is a fixed Markov family for  $\phi_t$
of small size $\chi < \ep_0/2 < 1$. Set  $U = \cup_{i=1}^k U_i$ and $\Intu (U) = \cup_{j=1}^k \Intu(U_j)$.
The {\it shift map} $\sigma : U   \longrightarrow U$ is given by $\sigma  = \piU \circ \pp$, where 
$\piU : R \longrightarrow U$ is the {\it projection} along stable leaves. 
Notice that  $\tau$ is constant on each stable leaf $W_{R_i}^s(x) = W^s_{\ep_0}(x) \cap R_i$. 
For any integer $m \geq 1$ and any function $h : U \longrightarrow \C$ define 
$h_m : U \longrightarrow \C$ by $h_m(u) = h(u) + h(\sigma(u)) + \ldots + h(\sigma^{m-1}(u))$.
 
Denote by $\widehat{U}$ the {\it core} of  $U$, i.e. the set of those $x\in U$ such that 
$\pp^m(x) \in \Int(R) = \cup_{i=1}^k \Int(R_i)$ 
for all $m \in \Z$. It is well-known (see \cite{kn:B}) that $\hU$ is a residual subset of $U$ and has full
measure with respect to any Gibbs measure on $U$.
Clearly in general $\tau$ is not continuous on $U$, however $\tau$ is {\it essentially Lipschitz} on $U$ 
in the sense that there exists a constant $L > 0$ such that if $x,y \in U_i \cap \sigma^{-1}(U_j)$ 
for some $i,j$, then $|\tau(x) - \tau(y)| \leq L\, d(x,y)$.
The same applies to $\sigma : U \longrightarrow U$.  
%Throughout we will mainly  work with the restrictions of $\tau$ and $\sigma$ to $\hU$. 
%Set $\hU_i = U_i \cap \hU$.

Let $B(\hU)$ be the {\it space of  bounded functions} $g : \hU \longrightarrow \C$ with its standard norm  
$\|g\|_0 = \sup_{x\in \hU} |g(x)|$. Given a function $g \in B(\hU)$, the  {\it Ruelle transfer operator } 
$L_g : B(\hU) \longrightarrow B(\hU)$ is defined by $\di (L_gh)(u) = \sum_{\sigma(v) = u} e^{g(v)} h(v)\;.$
If $g \in B(\hU)$ is essentially Lipschitz on $\hU$, then  $L_g$ preserves the space $\clip(\hU)$
of {\it Lipschitz functions} $h: \hU \longrightarrow \C$.

The hyperbolicity of the flow on $M$ and the additional assumption (in section 4 below) that
the local stable holonomy maps are uniformly Lipschitz implies the existence of
 constants $c_0 \in (0,1]$ and $\gamma_1 > \gamma > 1$ such that
\begin{equation}
 c_0 \gamma^m\; d (u_1,u_2) \leq 
d (\sigma^m(u_1), \sigma^m(u_2)) \leq \frac{\gamma_1^m}{c_0} d (u_1,u_2)
\end{equation}
whenever $\sigma^j(u_1)$ and $\sigma^j(u_2)$ belong to the same  $U_{i_j}$ 
for all $j = 0,1 \ldots,m$.

%Sect. 3 
\section{\sc Comparison of ball diameters}
\setcounter{equation}{0}

\subsection{\sc Lyapunov regularity}
Let $M$ be a $C^2$ compact Riemann manifold and $\phi_t$ an Anosov  $C^2$  flow on $M$.
Set $f = \phi_1$ and denote by $\ll$ the set of all {\it Lyapunov regular points} of $f$ (see \cite{kn:P1}
or section 2.1 in  \cite{kn:BP}). It is well-known that $\ll$ is dense in $M$ and has full measure with respect to any
$f$-invariant probability measure on $M$. Let 
$$\lambda_1  < \lambda_2 < \ldots < \lambda_s$$ 
be the {\it exponentials of the positive Lyapunov exponents} of $f$ over $\ll$ (so we have $1 < \lambda_1$).
{\bf Fix an arbitrary constant} $\alpha \in (0,1]$ such that 
\be
\lambda_j^\alpha < \lambda_{j+1} \quad , \quad 1 \leq j < s\;.
\ee

Next, set $\nu_0 = (1+\lambda_1)/2$, $\nu_1 = \lambda_1 + (\lambda_2 - \lambda_1)/3$ and 
$\nu_2 = \lambda_1 + 2(\lambda_2 - \lambda_1)/3$ , so that 
$$1 < \nu_0 < \lambda_1 < \nu_1 < \nu_2 < \lambda_2\;.$$
Then take $\mu > 0$ so small that 
\be
\nu_0 e^{8\mu} < \lambda_1  \quad , \quad \nu_2 e^{8\mu} < \lambda_{2}
\ee
(it then follows that $ \lambda_1 e^{8\mu} < \nu_1$ and  $\nu_1 e^{8\mu} < \nu_{2}$), and
\be
0 < \mu < \min \left\{ \frac{\alpha}{2(2+\alpha)}\, \ln \nu_0 \:,
\; \ln \frac{\lambda_1 + \nu_1}{2\lambda_1}\;,\; \ln \frac{2\lambda_2}{\lambda_2 + \nu_2} \right\}\;.
\ee
Set
\be
\gamma = \max \{ (\nu_0/\lambda_1)^{\alpha} \;, \; \nu_1/\nu_2 \} < 1\;,
\ee
and consider
%\be
$$1 < \nu_0 < \tmu_1 = \lambda_1 e^{-2\mu} < \mu_1 = \lambda_1 e^{-\mu} 
< \lambda_1 <  \lambda'_1 = \lambda_1 e^{\mu}  < \tla_1 = \lambda_1 e^{2\mu}   <
\nu_1 < \nu_2 < \mu_2 = \lambda_2 e^{-\mu} < \lambda_2\;.$$
%\ee

Fix for a moment $\mu > 0$ with the above properties. 
%(some additional conditions on  $\mu$ will be imposed later). 
Then, for $x\in \ll$ we have an $f$-invariant decomposition
$$E^u(x) = E^u_1(x) \oplus E^u_2(x) \oplus \ldots \oplus E^u_s(x)$$
into subspaces of constant dimensions $n_1, \ldots, n_s$ such that for some 
{\it Lyapunov $\mu$-regularity function}  $R = R_{\mu} : \ll \longrightarrow (1,\infty)$,
i.e. a function with
\be
e^{-\mu} \leq \frac{R (f(x))}{R (x)} \leq  e^{\mu} \quad , \quad x\in \ll\;,
\ee
we have
\be
\frac{1}{R (x)\, e^{n\mu}} \leq \frac{\|df^n(x)\cdot v\|}{\lambda_i^n\|v\|} 
\leq R (x)\, e^{n\mu} \quad , \quad x\in \ll 
\;, \; v\in E^u_i(x)\setminus \{0\} \;, \; n \geq 0\;.
\ee

For $x \in \ll$ set 
$$\tE^u_2(x) = E^u_2(x) \oplus \ldots \oplus E^u_s(x)\;.$$ 
For any $u\in E^u(x)$ we will write $u = \uo + \ut$, where $\uo \in E^u_1(x)$ and $\ut\in \tE^u_2(x)$.
We will denote by $\|\cdot\|$ the norm on $E^u(x)$ generated by the Riemann metric, and we will
also use the norm  $|u| = \max\{ \|\uo\| , \|\ut\|\}$.
Taking the regularity function $R(x)$ appropriately (see \cite{kn:P1},\cite{kn:BP} or \cite{kn:PS}), 
we may assume that  
$$|u| \leq \|u\| \leq R(x) |u| \quad , \quad x\in \ll\;,\; u\in E^u(x)\;.$$

It follows from the general theory of partial hyperbolicity (see \cite{kn:P1}, \cite{kn:P2}, \cite{kn:BP}) 
that the invariant bundle 
$\{\tE^u_2(x)\}_{x\in \ll}$ is uniquely integrable over $\ll$, i.e. 
there exists a continuous $f$-invariant family $\{ \Wut_{\tq(x)}(x)\}_{x\in \ll}$ 
of $C^2$ submanifolds
$\Wut_{\tq(x)}(x)$ of $M$ tangent to the bundle $\tE^u_2$ for some 
Lyapunov {\it $\mu/2$-regularity function} $\tq = \tq_{\mu/2} : \ll \longrightarrow (0,1)$.
Moreover, it follows from Theorem 6.6 in \cite{kn:PS} and (3.1) that there exists
an $f$-invariant family $\{ \Wuo_{\tq(x)}(x)\}_{x\in \ll}$ of
$C^{1+\alpha}$ submanifolds $\Wuo_{\tq(x)}(x)$ of $M$ tangent to the bundle $\tE^u_1$.
(However this family is not unique in general.) For each $x\in \ll$ fix
an $f$-invariant family $\{ \Wuo_{\tq(x)}(x)\}_{x\in \ll}$ with the latter properties.
Then we can find a Lyapunov $\mu$-regularity function  $q = q_\mu: \ll \longrightarrow (0,1)$ 
and for any $x\in \ll$ a $C^{1+\alpha}$ diffeomorphism
$$\Phi_x : E^u (x; q(x)) \longrightarrow \Phi_x (E^u (x; q (x)) \subset W^u_{\tq (x)}(x)$$
such that 
\be
\Phi_x(E^u_1(x;q(x))) \subset \Wuo_{\tq(x)}(x)\quad , \quad
\Phi_x(\tE^u_2(x;q(x))) \subset \Wut_{\tq (x)}(x) \quad x\in \ll\;.
\ee
We will assume without loss of generality that the regularity function $R$ satisfies
\be
\|d \Phi_x(u)\| \leq R(x) \quad ,\quad  \|(d\Phi_x(u))^{-1}\| \leq R(x) \quad, 
\quad x\in \ll \:, \: u \in E^u (x ; q(x))\;.
\ee
For any $x\in \ll$ consider the $C^{1+\alpha}$ map (defined locally near $0$)
$$\hf_x = (\Phi_{f(x)})^{-1} \circ f\circ \Phi_x : E^u (x)  \longrightarrow E^u (f(x))\;.$$

Given $y \in \ll$ and any integer $k \geq 1$ we will use the notation
$$\hf_y^k = \hf_{f^{k-1}(y)} \circ \ldots \circ \hf_{f(y)} \circ \hf_y\quad,
\quad \hf_y^{-k} = (\hf_{f^{-k}(y)})^{-1} \circ \ldots \circ (\hf_{f^{-2}(y)})^{-1} 
\circ (\hf_{f^{-1}(y)})^{-1} \;,$$
at any point where these sequences of maps are well-defined.

It is well known (see e.g. the Appendix in \cite{kn:LY} or section 3 in \cite{kn:PS})  
that there exist Lyapunov $\mu$-regularity functions 
$\Gamma = \Gamma_\mu : \ll \longrightarrow [1,\infty)$ and $q = q_\mu : \ll \longrightarrow (0,1)$ and 
for each $x\in \ll$ a norm $\| \cdot \|'_x$ on $T_xM$ such that
\be
\|v\| \leq \|v\|'_x \leq \Gamma (x) \|v\| \quad ,\quad x\in \ll\:,\: v \in T_xM\;,
\ee
and for any $x \in \ll$ and any integer $n \geq 0$, assuming 
$\hf_x^j(u), \hf_x^j(v) \in E^u(f^j(x),q (f^j(x)))$ are well-defined for all 
$j =1, \ldots,n$, the following hold:
\be
\mu^n_2\, \|u-v\|'_x  \leq \|\hf_x^n(u) - \hf_x^n(v)\|'_{f^n(x)} \quad , 
\quad u,v \in \tE^u_2(x;q(x))\;,
\ee
\be
\mu^n_1\, \|u-v\|'_x  \leq \|\hf_x^n(u) - \hf_x^n(v)\|'_{f^n(x)} \leq (\lambda'_1)^n\, \|u-v\|'_x
\quad , \quad u,v \in E^u_1(x;q(x))\;,
\ee
\be
\mu^n_1\, \|u-v\|'_x  \leq \|\hf_x^n(u) - \hf_x^n(v)\|'_{f^n(x)} 
\quad , \quad u,v \in E^u(x;q(x))\;,
\ee
\be 
\mu_1^n \, \|v\|'_x \leq \|d\hf_x^n(u)\cdot v\|'_{f^n(x)} \leq (\lambda'_1)^n \, \|v\|'_x  \quad , \quad x\in \ll \:,\:
 u \in E^u(x;q(x))\;,\; v \in E^u_1(x)\;,
\ee
and
\be 
\mu_2^n \, \|v\|'_x \leq \|d\hf_x^n(u)\cdot v\|'_{f^n(x)}   \quad , \quad x\in \ll \:,\:
 u \in E^u(x;q(x))\;,\; v \in \tE^u_2(x)\;.
\ee
Clearly each of the above inequalities provides a corresponding inequality involving the norm $\|\cdot \|$. For example
(3.9) and (3.10) imply
\be
\mu^n_2\, \|u-v\|  \leq \mu^n_2\, \|u-v\|'_x  \leq \|\hf_x^n(u) - \hf_x^n(v)\|'_{f^n(x)} 
\leq \Gamma (f^n(x)) \|\hf_x^n(u) - \hf_x^n(v)\|\;
\ee
for all $x\in \ll$ and all $u,v\in \tE^u_2(x)$.

\def\Bmt{\overline{B_{\ep_0}(\mt)}}

\subsection{\sc Balls in Bowen's metric}

We will use the notation from section 3.1. 
Given  $t > 0$ and  $\delta > 0$ set
$$B^u_t (x,\delta) = \{ y\in W^u_{\delta}(x) : d (\phi_t(x), \phi_t(y)) \leq \delta  \}\;.$$

%(???? Different balls??)

Our aim in this section is to prove the following

\bs

\noindent
{\bf Theorem 3.1.} {\it  For every $\ep > 0$ there exist Lyapunov $\ep$-regularity functions
$\omega : \ll \longrightarrow (0,1)$ and $G: \ll \longrightarrow [1,\infty)$ such that
for  any $0 < \delta_1 < \delta_2$ there exists a constant $K = K(\delta_1,\delta_2) \geq 1$
with the following property: for all  $x\in \ll$ with $\delta_2 \leq \omega(x)$
and all $ t > 0$ we have}
$$\diam(B^u_t( \phi_{-t}(x), \delta_2)) \leq K\, G(x) \diam (B^u_t( \phi_{-t}(x),\delta_1))\;.$$
%where $z = \phi_{-t}(x)$.}

\ms

{\bf Fix an arbitrary} $\ep > 0$. Let $\mu > 0$ satisfy (3.2), (3.3) and
$$ 0 < \mu <  \frac{\ep \alpha}{12}\;.$$

\def\chBo{\check{B}^{u,1}}
\def\tBo{\tB^{u,1}}
\def\hBo{\hB^{u,1}}
\def\hpi{\hat{\pi}}
\def\tU{\widetilde{U}}
\def\tr{\tilde{r}}

For a non-empty set $X\subset E^u (x)$  set
$$\ell(X) = \sup \{ \|u\|  : u \in X\}\;.$$

Given $z\in \ll$ and $p \geq 1$, setting $x = f^p(z)$, define
$$\hB^u_p (z, \delta) = \{ u \in E^u (z) : \|\hf^p_z(u) \| \leq\delta\}\quad, \quad
\hBo_p (z, \delta) = E^u_1(z) \cap \hB^u_p(z,\delta)\;,$$
%\{ u \in E^u_1 (z) : \|\hf^p_z(u) \| \leq\delta\}\;,$$
$$\tB^u_p (z, \delta) = \{ u \in E^u (z) : \|d \hf^p_z(0) \cdot u\| \leq\delta\}\quad ,\quad
\tBo_p (z, \delta) = E^u_1(z) \cap \tB^u_p(z,\delta)\;.$$
%\{ u \in E^u_1 (z) : \|d \hf^p_z(0) \cdot u\| \leq\delta\}\;.$$
Notice that $u \in E^u_1(z)$ implies $\hf_z^p(u) \in E^u_1(f^p(z))$ whenever
$\hf_z^p(u)$ is well-defined.

Theorem 3.1 will be  derived from Lemma 3.3 below and the following proposition.

\bs

\noindent
{\bf Proposition  3.2.} {\it There exists a $12\mu/\alpha$-regularity function 
$\omega : \ll \longrightarrow (0,1)$ with 
$\omega(x) \leq q(x)$ for all $x\in \ll$ and a $4\mu$-regularity function 
$G : \ll \longrightarrow [1,\infty)$ such that
for any  $x\in \ll$, any $\delta \in (0,  \omega(x)]$ and any integer $p \geq 1$ for  $z = f^{-p}(x)$ we have}
$$\ell(\hB^u_p (z, \delta)) \leq G(x) \ell(\hBo_p (z, \delta))\;.$$

\ms

The proof of Proposition 3.2 takes most of this section. 

Taylor's formula (see also section 3 in \cite{kn:PS}) implies that there exists 
a Lyapunov $\mu$-regularity function  $D = D_\mu : \ll \longrightarrow  [1,\infty)$ such that 
for any $i = \pm 1$ we have
\be
\quad\|\hf^{i}_x(v) - \hf^{i}_x(u) - d\hf^{i}_x(u) \cdot (v-u)\| \leq D(x) \, \|v-u\|^{1+\alpha} \:\:,\: 
x\in \ll\:, \: u,v \in E^u (x ; q(x))\;.
\ee
%and
%\be
%\|d\hf_x(u) - d\hf_x(0)\| \leq D_\mu(x)\, (\|u\|)^\alpha \quad , \quad x\in \ll \:, \: u \in E^u (x;q_\mu(x))\;.
%\ee

Fix for a moment $x\in \ll$ and an integer $p \geq 1$, set $z = f^{-p}(x)$ and given
$v \in E^u (z; q(z))$, set
\be 
z_j = f^j(z)\quad, \quad v_j = \hf^j_z(v) \in E^u (z_j) \quad , \quad w_j = d\hf_z^j(0)\cdot v \in E^u (z_j)
\ee
for any $j = 0,1,\ldots, p$ (assuming that these points are well-defined).

%The following is the main step in the proof of Theorem 3.1. Its proof uses an argument
%from the proof of Lemma 3.3 in \cite{kn:St3}.

\bs

\noindent
{\bf Lemma 3.3.} {\it There exist a Lyapunov $6\mu$-regularity function
$L = L_{6\mu}: \ll \longrightarrow [1,\infty)$ and a Lyapunov $7\mu/\alpha$-regularity 
function $\tr = \tr_{7\mu/\alpha}: \ll \longrightarrow (0,1)$ with $\tr \leq q$
%and a global constant $D' \geq 1$ 
such that for any $x\in \ll$, any integer $p \geq 1$ and any $v\in E^u(z,\tr (z))$ with 
$\|\hf_z^p(v)\|\leq \tr(x)$,  where $z = f^{-p}(x)$, we have
$$\|\wo_p - \vo_p\| \leq L(x) |v_p|^{1+\alpha}\;.$$
Moreover, if $|v_p| = \|\vo_p\| \neq 0$, then}  $1/2 \leq \|\wo_p\|/\|\vo_p\| \leq 2$.

\bs

\noindent
{\it Proof of Lemma} 3.3. One checks easily that
\be
\tr(x) = \left(\frac{1-\gamma}{2}\right)^{1/\alpha} \cdot \frac{q(x)}{D(x)^{1/\alpha}
\Gamma(x)^{2/\alpha + 1} R(x)^{1/\alpha+1}} \leq q(x)
\ee
defines a Lyapunov $7\mu/\alpha$-regularity function on $\ll$.

Let $x\in \ll$ and $z = f^{-p}(x)$ for some integer $p \geq 1$, and let $v\in E^u(z,\tr (z))$  
be such that $\|\hf_z^p(v)\|\leq \tr(x)$. Using the notation (3.17), by (3.9) and (3.12),
%\be
$$\|v_k\| = \|\hf^k_z(v)\|  \leq \frac{\Gamma(x)}{\mu_1^{p-k}}\, \|\hf_z^p(v)\|
= \frac{\Gamma(x)}{\mu_1^{p-k}}\, \|v_p\|\;$$
%\ee
for all $ k = 0,1, \ldots,p-1$. It follows from (3.16) that 
$\|\hf_z(v) - d\hf_z(0)\cdot v\| \leq D \, \|v\|^{1+\alpha}$, so
$w_1  = \hf_z (v) + u_1 = v_1+u_1$
for some $u_1 \in E^u (z_1)$ with $\|u_1\|\leq D(z_1)\, \|v\|^{1+\alpha}$. Hence
$w_2 = d\hf_{z_1}(0)\cdot w_1 = d\hf_{z_1}(0)\cdot v_1 + d\hf_{z_1}(0)\cdot u_1$.
Using (3.16) again, we get $d\hf_{z_1}(0)\cdot v_1 = \hf_{z_1}(v_1) + u_2 = v_2 + u_2$
for some $u_2 \in  E^u (z_2)$ with $\|u_2\|\leq D(z_2)\, \|v_1\|^{1+\alpha}$. 
Thus, $w_2 = v_2 + u_2 + d\hf_{z_1}(0)\cdot u_1$.
Continuing by induction, as in the proof of Lemma 3.3 in \cite{kn:St3}, one derives
\be
w_p = v_p + u_p + d\hf_{z_{p-1}}(0)\cdot u_{p-1} + d\hf^2_{z_{p-2}}(0)\cdot u_{p-2} +
\ldots + d\hf_z^{p-1}(0)\cdot u_1\;,
\ee
where $u_j \in E^u (z_j)$ and $\|u_j\| \leq D (z_j) \, \|v_{j-1}\|^{1+\alpha}$ for all $j = 1,\ldots,p$.
Then 
$$\|u_j\| \leq D (z_j)\, \|v_{j-1} \|^{1+\alpha} \leq  
\frac{D(x) \, e^{(p-j)\mu} \Gamma(x)^{1+\alpha}}{\mu_1^{(1+\alpha)(p-j)}}\, \|v_p\|^{1+\alpha}
\leq \frac{D(x) \, \Gamma (x)^{1+\alpha}}{\tmu_1^{(1+\alpha)(p-j)}}\, \|v_p\|^{1+\alpha}\;.$$
Combining the latter with (3.9) and (3.13) gives
\begin{eqnarray*}
\|d\hf_{z_j}^{p-j}(0)\cdot \uo_j\|
 \leq  (\lambda'_1)^{p-j} \, \Gamma(z_j) \|u_j\| 
 \leq  D (x)\, \Gamma(x)^{2+\alpha}\; \left(\frac{\tla_1}{\tmu_1^{1+\alpha}}\right)^{p-j}\, 
 \|v_p\|^{1+\alpha}\;
 \end{eqnarray*}
for all $j = 1, \ldots, p$. Setting $L'(x) = D (x) \Gamma(x)^{2+\alpha}$, and using (3.3) and (3.4)
to get
$$\frac{\tla_1}{\tmu_1^{1+\alpha}} = \frac{\lambda_1 e^{2\mu}}{\lambda_1^{1+\alpha} e^{-2\mu(1+\alpha)}}
= \frac{e^{2\mu(2+\alpha)}}{\lambda_1^\alpha} \leq \left(\frac{\nu_0}{\lambda_1}\right)^\alpha 
\leq  \gamma\;,$$
it follows that $\|d\hf_{z_j}^{p-j}(0)\cdot \uo_j\|\leq  L'(x)\, \gamma^{p-j} \, \|v_p\|^{1+\alpha}$. Now
(3.19) yields
$$\|\wo_p - \vo_p\| 
\leq  L' (x) \, \| v_p\|^{1+\alpha} \, \sum_{j=1}^p \gamma^{p-j}
\leq \frac{L'(x)}{1-\gamma} \, \| v_p\|^{1+\alpha} 
\leq \frac{L'(x) R(x)^{1+\alpha}}{1-\gamma} \, | v_p|^{1+\alpha} \;.$$
Since $L(x) = L'(x)R(x)^{1+\alpha}/(1-\gamma)$ is a Lyapunov $6\mu$-regularity function, this proves
the first part of the lemma.

If $|v_p| = \|\vo_p\|$, then the above gives $\|\wo_p - \vo_p\| \leq L(x) \|\vo_p\|^{1+\alpha}$, so
$$\left| \frac{\|\wo_p\|}{\|\vo_p\|} - 1\right| \leq L(x) \|\vo_p\|^\alpha \leq 
L(x) (\tr(x))^{\alpha} \leq \frac{1}{2}\;,$$
by the choice of $\tr(x)$. Hence $1/2 \leq \|\wo_p\|/\|\vo_p\| \leq 2$.
\endofproof

\bs

\noindent
{\bf Corollary 3.4.} {\it Under the assumptions of Lemma 3.3, for $u = \vo \in E^u_1(z)$ we have}\\
$1/2 \leq \|\wo_p\|/\|u_p\| \leq 2$.

\bs

\noindent
{\it Proof.} We just apply Lemma 3.3 replacing $v$ by $u$. Since $\hf^p_z(u) \in E^u_1(x)$, we have
$u_p = \uo_p$.
\endofproof

\bs

\def\tde{\tilde{\delta}}

\noindent
{\bf Lemma 3.5.} {\it Assume that the regularity function $r$ satisfies $r(x) \leq \tr(x)$ and
\be
r(x) \leq  \min \left\{ 
\left( \frac{1/\nu_2 - 1/\lambda_2}{2 e^\mu \Gamma^2 (x) D (x)} \right)^{1/\alpha}  \;, \;
\left( \frac{1/\lambda_1 - 1/\nu_1}{2 e^{3\mu} \Gamma^2(x) D(x)} \right)^{1/\alpha} \right\}\;
\ee
for all $x\in \ll$.
Then for any $x\in \ll$ and any $V = \Vo + \Vt \in E^u(x; r(x))$ we have
\be
\|(\hf_x^{-1})^{(2)}(V)\|'_{f^{-1}(x)} \leq \frac{\|\Vt\|'_x}{\nu_2} \;,
\ee
and}
\be
\|(\hf_x^{-1})^{(1)}(V)\|'_{f^{-1}(x)} \geq \frac{\|\Vo\|'_x}{\nu_1} \;.
\ee

\noindent
{\it Proof of Lemma} 3.5. Let $x\in \ll$, $y = f^{-1}(x)$ and let 
$V = \Vo + \Vt \in E^u(x; r(x))$. 
%be such that $\|V\|'_x \leq \ep_\mu$. 
By (3.16),
$$\hf_x^{-1}(V) - \hf_x^{-1}(\Vo,0) = d\hf_x^{-1}(\Vo,0)\cdot (0,\Vt) + \xi$$
for some $\xi\in E^u(y)$ with  $\|\xi\| \leq D(y) \|\Vt\|^{1+\alpha}$. 
Using $\|\Vt\| \leq r(x)$ and (3.9), the latter gives
$$\|\xi\|'_y \leq \Gamma(y) \|\xi\| \leq \Gamma(y) D(y)  \|\Vt\|^{1+\alpha} \leq 
\Gamma(x) D(x)e^{2\mu} \|\Vt\|^{1+\alpha}
\leq \Gamma^2(x) D(x) e^{2\mu} \|\Vt\|'_x \, r^\alpha(x)\;.$$
Since
%$$\hf_x^{-1}(\Vo,0) = d\hf_x^{-1}(0)\cdot (V_1,0) + w'\;,$$
$\hf_x^{-1}(\Vo,0) \in E^u_1(y)$, it follows from (3.10) that
\begin{eqnarray}
\|(\hf_x^{-1})^{(2)}(V)\|'_y 
& \leq & \|d\hf_x^{-1}(\Vo,0)\cdot (0,\Vt)\|'_y + \|\xi\|'_y\nonumber\\
& \leq & \|\Vt\|'_x \; \left(\frac{1}{\mu_2} + \Gamma^2(x) D(x) e^{2\mu}  r^\alpha(x)\right)\;.
\end{eqnarray}
Now (3.9) and (3.20) imply
\begin{eqnarray*}
\frac{1}{\mu_2} + \Gamma^2(x) D(x) e^{2\mu}  r^\alpha(x)
 \leq  \frac{e^\mu}{\lambda_2} + e^\mu\, \frac{1/\nu_2 - 1/\lambda_2}{2}
= e^\mu \frac{\lambda_2 + \nu_2}{2 \lambda_2 \nu_2} < \frac{1}{\nu_2}\;,
\end{eqnarray*}
since by (3.3) we have $e^\mu < \frac{2\lambda_2}{\lambda_2 + \nu_2}$. 
The above and (3.23) imply (3.21).

Similarly, we have
$$\hf_x^{-1}(V) - \hf_x^{-1}(0,\Vt) = d\hf_x^{-1}(0,\Vt)\cdot (\Vo,0) + \eta$$
for some $\eta\in E^u(y)$ with $\|\eta\| \leq D(y) \|\Vo\|^{1+\alpha}$.
Then $\|\eta\|'_y \leq \Gamma^2(x) D(x) e^{2\mu} \|\Vo\|'_x \, r^\alpha(x)$. 
Since $\hf_x^{-1}(0,\Vt) \in \tE^u_2(y)$, by (3.11),
\begin{eqnarray}
\|(\hf_x^{-1})^{(1)}(V)\|'_y 
& \geq & \|d\hf_x^{-1}(0,\Vt)\cdot (\Vo,0)\|'_y - \|\eta\|'_y \nonumber\\
& \geq &  \|\Vo\|'_x \; \left(\frac{1}{\lambda'_1} -  \Gamma^2(x) D(x) e^{2\mu}  r^\alpha(x)\right) 
\;.
\end{eqnarray}
Now (3.20) implies
\begin{eqnarray*}
\frac{1}{\lambda'_1} -   \Gamma^2(x) D(x) e^{2\mu}  r^\alpha(x)
 \geq  \frac{1}{\lambda_1 e^\mu } - \frac{1/\lambda_1 - 1/\nu_1}{2e^\mu}
= \frac{\lambda_1 + \nu_1}{2 e^\mu \lambda_1 \nu_1} > \frac{1}{\nu_1}\;,
\end{eqnarray*}
since by (3.3) we have $e^\mu < \frac{\lambda_1 + \nu_1}{2\lambda_1}$. 
The above and (3.24) imply (3.22).
\endofproof

\bs

\noindent
{\it Proof of Proposition 3.2.} Define $\omega : \ll \longrightarrow (0,1)$ by
\be
\omega(x) =  \left(\frac{1-\gamma}{2}\right)^{1/\alpha}  
\frac{q(x)}{16 \nu_1 \Gamma(x)^{2/\alpha + 3} R(x)^{1/\alpha +1}}
\left( \frac{1/\nu_2 - 1/\lambda_2}{2 e^{3\mu} D(x)} \right)^{1/\alpha} 
\quad \;.
\ee
Clearly $\omega(x) \leq \tr(x)$, the function defined by (3.18). Moreover,
$1/\lambda_1 - 1/\nu_1 > 1/\nu_2 - 1/\lambda_2$ shows that 
$r(x) = 16\nu_1 \Gamma^3(x) R(x)\omega(x)$ satisfies
(3.20), so Lemma 3.5 applies. It is easy to check that $\omega$ is a 
Lyapunov $12\mu/\alpha$-regularity function.

Let $x\in \ll$ and let $p \geq 1$ be an integer. Set $z = f^{-p}(x)$. Given  $\delta > 0$
with $\delta \leq \omega(x)$, we have $16\nu_1 \Gamma^3(x) \delta \leq r(x)$. 
Let $v\in \hB^u_p(z,\delta)$ be such that $\|v\|$ is the maximal possible. 
Then we must have $\|v_p\| = \delta$,
where we use the notation in (3.17).  Set $V = v_p$.
%Notice that if
%$\|(\hf_x^{-j})^{(2)}(V)\|'_{x_{j}} < \|(\hf_x^{-j})^{(1)}(V)\|'_{x_j}$
%for some $j = 0,1,\ldots,p-1$,  then it follows from (3.21) and (3.22) that
%$\|(\hf_x^{-i})^{(2)}(V)\|'_{x_i} < \|(\hf_x^{-i})^{(1)}(V)\|'_{x_i}$ for all $i$ with
%$j \leq i \leq p$. 

Next, consider two cases.

\ms

{\bf Case 1.} $\|\vt\| \geq \|\vo\|$.
Let $U = (\Uo,0)$ be such that $\Uo \in E^u_1(x)$ is
an arbitrary element with $\|\Uo\| = \delta$. Then $u = \hf_x^{-p}(U)\in \hBo_p(z,\delta)$,
and applying (3.22) $p$ times gives
$\|u\|'_z \geq \frac{\|\Uo\|'_x}{\nu_1^p} \geq \frac{\delta}{\nu_1^p}$.
Similarly, applying (3.21) $p$ times and using (3.4) implies
$$\|\vt\|'_z \leq \frac{\|\Vt\|'_x}{\nu_2^p} \leq \frac{\Gamma(x) \|\Vt\|}{\nu_2^p}
\leq \Gamma(x) \frac{\delta}{\nu_2^p}
\leq \Gamma(x) \gamma^p\, \frac{\delta}{\nu_1^p} \leq \Gamma(x) \gamma^p\, \|u\|'_z\;.$$
This and (3.9) give
$$\|\vt\| \leq \|\vt\|'_z \leq \Gamma(x)\gamma^p\, \|u\|'_z
\leq \Gamma(x) \gamma^p \Gamma(z) \|u\|\leq \Gamma(x) \gamma^p e^{p\mu} \Gamma(x) \|u\|
\leq \Gamma^2(x) \|u\|\;.$$
The latter yields $\|\vo\| \leq \|\vt\| \leq \|\vt\|'_z \leq \Gamma^2(x) \|u\|$, and therefore
$|v| \leq \Gamma^2(x) \|u\|$. Hence
$\|v\| \leq R(x) |v| \leq \Gamma^2(x) R(x) \|u\|$, which shows that
$\ell(\hB^u_p(z,\delta)) \leq \Gamma^2(x)R(x) \ell(\hBo_p(z,\delta))$.

\ms

{\bf Case 2.} $\|\vt\| < \|\vo\|$. Set $x_j = z_{p-j} = f^{-j}(x)$. 
Let $q < p$ be the largest integer  with 
$\|(\hf_x^{-q})^{(2)}(V)\| \geq \|(\hf_x^{-q})^{(1)}(V)\|$.
Set $\tp = p-q$, $\tx = x_{\tp}$, 
%$\tde = \|\hf_z^{\tp}(v)\|'_{\tx}$, 
$u = (\vo,0) \in E^u_1(z)$ and $u_j = \hf_z^j(u)$. Then clearly $\|u\| = \|\vo\|$, 
and by the choice of $q$, $\|\vt_{\tp}\| \leq \|\vo_{\tp}\|$. Thus, $|v_{\tp}| = \|\vo_{\tp}\|$, 
and now it follows from Lemma 3.3 and Corollary 3.4 with $p$ replaced by $\tp$ that 
$\|u_{\tp}\| \leq 4 \|\vo_{\tp}\|$, so 
$\|u_{\tp}\|'_{\tx} \leq 4 \Gamma(\tx) \|\vo_{\tp}\|'_{\tx}$. Again by the choice of $q$,
$\|\vt_{\tp+1}\| \geq \|\vo_{\tp+1}\|$, so
$$\|\vo_{\tp+1}\|'_{x_{\tp+1}} \leq \Gamma(x_{\tp+1}) \|\vo_{\tp+1}\|
\leq \Gamma(x_{\tp+1}) \|\vt_{\tp+1}\| \leq \Gamma(x_{\tp+1}) \|\vt_{\tp+1}\|'_{x_{\tp+1}}
\leq \Gamma(x) e^{q\mu} \|\vt_{\tp+1}\|'_{x_{\tp+1}}\;.$$
Hence, using (3.12), (3.11), (3.21) and (3.9), we get
\begin{eqnarray*}
\|u_{\tp+1}\|'_{x_{\tp+1}} 
& \leq & \lambda'_1\|u_{\tp}\|'_{\tx} \leq 4 \lambda'_1 \Gamma(\tx) \|\vo_{\tp}\|'_{\tx}
\leq \frac{4 \lambda'_1}{\mu_1} \Gamma(\tx)\|\vo_{\tp+1}\|'_{x_{\tp+1}}
\leq 4  \lambda'_1 \Gamma^2(x) e^{2q\mu} \|\vt_{\tp+1}\|'_{x_{\tp+1}}\\
& \leq & \frac{4\lambda'_1}{ \nu_2^{q-1}} \Gamma^2(x) e^{2q\mu} \|\Vt\|'_x
 \leq  \frac{4\lambda'_1\, \Gamma^3(x) e^{2q\mu} \delta}{\nu_2^{q-1}}\;.
\end{eqnarray*}
Now (3.11) and  (3.4) imply
$$\|u_p\|'_x \leq (\lambda'_1)^{q-1}\, \|u_{\tp+1}\|'_{x_{\tp+1}}
\leq  (\lambda'_1)^{q-1}\,\frac{4\lambda'_1\, \Gamma^3(x) e^{2q\mu} \delta}{\nu_2^{q-1}}
\leq 4  \nu_1\, \Gamma^3(x)\delta\, \gamma^{q-1}
\leq 4 \nu_1 \Gamma^3(x)\, \delta\;,$$
and by (3.9), $\|u_p\| \leq 4\nu_1 \Gamma^3(x)\, \delta$.
Thus, $u \in \hBo_p(z,4\nu_1 \Gamma^3(x)\delta)$. 

Since $\|u\| = \|\vo\| > \|\vt\|$, it follows that $\|u\| \geq |v|$, so
$$\ell(\hB^u_p(z,\delta)) = \|v\| \leq R(x)|v| \leq R(x)\|u\| 
\leq R(x) \ell( \hBo_p(z,4\nu_1 \Gamma^3(x)\delta))\;.$$

On the other hand, it follows from Lemma 3.3, Corollary 3.4 and the linearity
of the map $d\hf_z^p(0)$ that 
$\ell(\hBo_p(z,4\nu_1 \Gamma^3(x)\delta)) \leq 16 \nu_1 \Gamma^3(x)\, 
\ell(\hBo_p(z,\delta))$. Thus, 
$\ell(\hB^u_p(z,\delta)) \leq  16 \nu_1 \Gamma^3(x) R(x)\,  \ell(\hBo_p(z,\delta))$.

It follows from cases 1 and 2 that the regularity functions $\omega(x)$ and
$G(x) = 16 \nu_1 \Gamma^3(x)R(x)$ satisfy the requirements of the proposition.
\endofproof

\bs

\noindent
{\it Proof of Theorem} 3.1. Clearly it is enough to prove the analogous statement
for sets of the form $\hB^u_p(z,\delta)$, and instead of diameters it is enough to
work with $\ell(\cdot)$.

Assume that the regularity functions $\omega$ and $G(x)$ are as in Proposition 3.2.

Let $x\in \ll$, $p \geq 1$ be an integer and $z = f^{-p}(x)$. We will use again the
notation (3.17). Let $0 < \delta_1 < \delta_2 \leq \omega(x)$.
%where $\ep_\mu$ is defined by (3.18). 
It follows from Proposition 3.2
that $\ell(\hB^u_p(z,\delta_2)) \leq G(x) \ell(\hBo_p(z, \delta_2))$, while the second part
of Lemma 3.3 shows that $\ell(\hBo_p(z,\delta_2)) \leq \ell(\tBo_p(z,2\delta_2))$.
Next, $\ell(\tBo_p(z,2\delta_2)) = \frac{4\delta_2}{\delta_1}\, \ell(\tBo_p(z,\delta_1/2))$
by the linearity of the map $d\hf_z^p(0)$.
Then using again the second part of Lemma 3.3, we get
$\ell(\tBo_p(z,\delta_1/2)) \leq \ell(\hBo_p(z,\delta_1)) \leq \ell(\hB^u_p(z,\delta_1))$.
Combining all these inequalities gives 
$\ell(\hB^u_p(z,\delta_2)) \leq K\, G(x) \ell(\hB^u_p(z,\delta_1))$, where
$K = \frac{4\delta_2}{\delta_1}$.
\endofproof

\bs

\noindent
{\it Proof of Theorem} 1.1. Consider the Anosov flow $\psi_t = \phi_{-t}$ on $M$. Clearly this flow
has the same set $\ll$ of Lyapunov regular points. Let $\omega$ and $G$ be Lyapunov regular
functions satisfying the requirements of Theorem 3.1 for the flow $\psi_t$.
We will denote by $\W^s_\delta(x)$ and 
$\W^u_\delta(x)$ the local stable and unstable manifolds for the flow $\psi_t$. Clearly
$W^s_\delta(x) = \W^u_\delta(x)$ and $W^u_\delta(x) = \W^s_\delta(x)$.

Given $x\in \ll$, $t > 0$ and $0 < \delta_1 < \delta_2  \leq \omega(x)$, set $x' = \phi_t(x)$, and notice that
$\phi_t(B^s(x,\delta_i)) = \B^u_t(x',\delta_i)$ for $ i = 1,2$, where
$$\B^u_t(x',\delta_i) = \{ y'\in \W^u_{\delta_i}(x') : d(\psi_t(x'), \psi_t(y')) \leq \delta_i\}\;.$$
Using  Theorem 3.1 for the flow $\psi$, it follows that there exists a constant $K(\delta_1, \delta_2) \geq 1$
such that $\diam(\B^u_t(x', \delta_2)) \leq KG(x) \diam (\B^u_t(x',\delta_1))$, i.e.
$\diam(\phi_t(B^s(x, \delta_2))) \leq KG(x) \diam (\phi_t(B^s(x,\delta_1)))$.
\endofproof

\def\Intu{\Int^u}
\def\Ints{\Int^s}
\def\tcc{\widetilde{\cc}}

%Sect. 4
\section{\sc Decay of cylinder diameters in a Markov coding}
\setcounter{equation}{0}

Let again $\rr = \{ R_i\}_{i=1}^k$ be a fixed Markov family as in section 2. Define the matrix 
$A = (A_{ij})_{i,j=1}^k$  by  $A_{ij} = 1$ if $\pp(\Int(R_i)) \cap \Int(R_j) \neq  \e$ and 
$A_{ij} = 0$ otherwise.  
According to \cite{kn:BR} (see section 2 there), we may assume that $\rr$ is chosen in such a way that 
$A^{M_0} > 0$ (all entries of the $M_0$-fold product of $A$  by itself are positive) for some  integer 
$M_0 > 0$. In what follows  we assume that the matrix $A$ has this property.

Given a finite string $\ii = (i_0,i_1, \ldots,i_m)$  of integers $i_j \in \{ 1, \ldots,k\}$, we will say 
that $\ii$ is {\it admissible}
if  for any $j = 0,1, \ldots,m-1$ we have $A_{i_j i_{j+1}} = 1$.  Given an admissible string $\ii$, 
denote by  $\co[\ii]$ the set
of those $x\in U$ so that $\sigma^j(x) \in \Intu(U_{i_j})$ for all $j = 0,1, \ldots,m$. The set 
$\di C[\ii] = \overline{\co[\ii]} \subset \mt$
will be called a {\it cylinder} of length $m$ in $U$, while $\co[\ii]$ will be called an {\it open cylinder} 
of length $m$. It follows from the properties of the Markov family that $\co[\ii]$ is an open dense 
subset of $C[\ii]$. Any cylinder of the form $C[i_0, i_1, \ldots, i_m, i_{m+1}, \ldots, i_{m+q}]$ 
will be called a {\it subcylinder} of $C[\ii]$ of {\it co-length} $q$.

In what follows the cylinders considered are always defined by finite admissible strings. 
Given $x\in U_i$ for some $i$ and $r > 0$ we will denote by $B_U(x,r)$ the set of all $y\in U_i$ with 
$d (x,y) < r$.

It is easy to see that $\diam(C[\ii]) \to 0$ exponentially fast as $m \to \infty$. 
A much more subtle question is if there exists a constant $\rho \in (0,1)$ such that 
for any cylinder $\cc = C[i_0,i_1, \ldots,i_m]$ and any subcylinder
$\cc' = C[i_0,i_1, \ldots,i_m, i_{m+1}]$ we have $\diam (\cc') \geq \rho \, \diam(\cc)$. 
Using Theorem 1.1 here we show that this is always the case under some regularity assumptions about the flow.

Recall the constants $c_0 \in (0,1)$ and  $\gamma_1 > \gamma > 1$ from section 2,  
and fix an integer  $p_1 \geq 1$ with
\be
\rho_0 = \frac{1}{c_0\gamma^{p_1}} < \min \left\{ \frac{\diam(U_i)}{\diam (U_j)}   
: i,j = 1, \ldots, k \right\}\;.
\ee
Then clearly $\rho_0 < 1$. Set $\rho_1 = \rho_0^{1/p_1}$
and   fix a constant $r_0 > 0$ with $2r_0 < \min \{ \diam(U_i) : i = 1, \ldots,k\}$ and for each 
$i = 1, \ldots,k$ a  point $\hz_i \in \hU_i$ such that $B_U (\hz_i, r_0) \subset \Intu(U_i)$.

The following is an easy consequence of (2.1).

\bs

\noindent
{\bf Lemma 4.1.}(\cite{kn:St1}) {\it There exists a global constant $C_1 > 0$ such that for any cylinder 
$C[\ii]$ of length $m$ we have
$\diam (C[\ii] ) \leq C_1\, \rho_1^m$ and $\diam (C[\ii] ) \geq \frac{c_0 r_0}{\gamma_1^m}$.}
\endofproof

\bs

In what follows we will assume that $\rho_1 \in (0,1)$ and $C_1 > 0$ are fixed constants
with the above property. Fix a constant $\ep > 0$  such that
\be
e^{-\ep} > \sqrt{\rho_1}\;.
\ee

From now on we will assume that the local stable holonomy maps through $\mt$ are uniformly Lipschitz.
Then there exists a constant $L \geq 1$ such that $d(\pi_y(z), \pi_y(z')) \leq L\, d(z,z')$ for all 
$x,y\in M$ with $d(x,y) < \ep_1$ and $z,z'\in W^u_{\ep_1}(x)$. (See section 2 for the choice
of $\ep_1$.)

Given $i = 1, \ldots, k$, according to the choice of the Markov family $\{R_i\}$,
the projection\\ $\pr_{R_i} : W_i = \phi_{[-\chi, \chi]}(R_i) \longrightarrow R_i$ along the flow $\phi_t$
is well-defined and Lipschitz.  Since the  projection $\pi_i : R_i \longrightarrow U_i$ along stable leaves 
is Lipschitz, the map $\psi_i = \pi_i\circ \pr_{R_i} : W_i  \longrightarrow W^u_{R_i}(z_i)$ is also 
Lipschitz. Thus, we may assume the constant $L \geq 1$ is chosen sufficiently large so that 
$d (\psi_i(u),\psi_i (v)) \leq L\, d (u,v)$ for all $u,v\in W_i$
and all $i = 1, \ldots,k$.

Next, if $V = W^u_R(x)$ is the unstable leaf of some point $x\in R_i$ and
$\ii = (i_0 = i,i_1, \ldots,i_m)$ is an admissible sequence, consider the {\it generalized
cylinder}
$$C_V[\ii] = \{ y\in V : \pp^j(y) \in R_{i_j} \:, \: j = 0,1,\ldots,m\}\;.$$
Clearly, $\pi_i(C_V[\ii]) = C[\ii]$, so
\be
\frac{1}{L} \diam (C_V[\ii]) \leq \diam(C[\ii]) \leq L\, \diam(C_V[\ii])\;
\ee
for any choice of $V$ and the admissible sequence $\ii$. For $V$ as above, $x \in V$ and $\delta > 0$ set
$$B_V(x,\delta) = \{ y \in V : d(x,y) < \delta\}\;.$$

\ms

\noindent
{\bf Theorem 4.2.} {\it Assume that $\phi_t : M \longrightarrow M$ is a $C^2$ Anosov flow 
such that the local stable holonomy maps are uniformly Lipschitz. Then there exist a 
constant $\rho \in (0,1)$  and a positive integer $p_0 \geq 1$  such that:}

(a) {\it  For  any cylinder $C[\ii] = C[i_0, \ldots,i_m]$ and any subcylinder
$C[\ii'] = C[i_0,i_{1}, \ldots, i_{m+1}]$ of $C[\ii]$ of co-length $1$  we have}
$\rho \; \diam ( C [\ii] ) \leq  \diam ( C [\ii'] ) \;.$

\ms

(b) {\it For  any cylinder $C[\ii] = C[i_0, \ldots,i_m]$ and any subcylinder
$C[\ii'] = C[i_0,i_{1}, \ldots, i_{m+1}, \ldots, i_{m+p_0}]$ of $C[\ii]$ of co-length $p_0$
we have} $\diam(C[\ii'] ) \leq \rho\, \diam (C[\ii]) \;.$

\medskip

\noindent
{\it Proof of Theorem} 4.2. Notice that  for any (admissible) $\ii$ we have 
$\sigma^m (\hC[i_0, \ldots,i_m]) = \hU_{i_{m}}$.

As in section 2 we will assume that the point $z_i \in \Int(R_i)$ is Lyapunov regular. 
Given $\ep > 0$ with (4.2), let  $\omega: \ll \longrightarrow (0,1)$ 
and $G : \ll \longrightarrow [1,\infty)$ be Lyapunov $\ep$-regularity functions 
with the properties described in Theorem 3.1.  
Fix a constant $r > 0$ with  $r \leq \min_{1\leq i \leq k} \omega(z_i)$ 
such that $B_U (z_i, r) \subset \Intu(U_i)$ for all $i = 1, \ldots,k$, and 
set $G_0 = \max_{1\leq i \leq k} G(z_i)$. Then fix an integer $p \geq p_1$ so large that
\be
\rho_1^{(p+1)/2} < \frac{r}{C_1 L}\;.
\ee

First note the following. Let $x\in R_j$ be such that $\pp^{p+1}(x) = z_i$ for some $j$ and $i$. Then
$z\in \cc = C_V[\ii]$ for some $\ii = [i_0, \ldots,i_{p+1}]$ with $i_0 = j$ and $i_{p+1} = i$, where
$V = W^u_R(z)$.  We claim that
\be
C_V[\ii] \subset B_V(z,\omega(x)) \subset V\;.
\ee
Indeed, by Lemma 4.1, $\diam(\cc) \leq C_1 \rho_1^{p+1}$. On the other hand, $\omega$ is a
Lyapunov  $\ep$-regularity function, so using (4.2), Lemma 4.1 and (4.3), we get
$$\omega(x) \geq \omega(z_i) e^{-(p+1)\ep} \geq r \, \rho_1^{(p+1)/2} 
= \frac{r}{\rho_1^{(p+1)/2}} \; \rho_1^{p+1}
> C_1L \, \rho_1^{p+1} \geq  L\, \diam(C[\ii]) \geq  \diam(C_V[\ii])\;.$$
This proves (4.5).

\medskip

(a) Assume that $m > p$, and let $\ii = [i_0,i_1, \ldots,i_m]$ and
$\ii' = [i_0,i_1, \ldots,i_m, i_{m+1}]$ be admissible sequences.
%Set $z = z_{i_{m+1}}$ and 
Let $z \in R_{i_0}$ be such that $\pp^{m+1}(z) = z_{i_{m+1}}$
and $\pp^j(z) \in R_{i_j}$ for all $j = 0,1,\ldots,m+1$. Set $V = W^u_R(z)$,
$\cc = C_V[\ii]$ and $\cc' = C_V[\ii']$. 

Next, set  $x = \pp^{m-p}(z)$ and $V' = W^u_R(x)$, and consider the cylinders
$\tcc = C_{V'}[i_{m-p}, i_{m-p+1}, \ldots, i_m]$ and   
$\tcc' = C_{V'}[i_{m-p}, i_{m-p+1}, \ldots, i_m, i_{m+1}]$. Since $\pp^{p+1}(x) = z_{i_{m+1}}$,
using (4.5) we get $\tcc \subset B_{V'}(x,\omega(x))$, and moreover $\tcc \subset B_{V'}(x,C_1\rho_1^{p+1})$.
On the other hand it is easy to see using (2.1) that $\tcc' \supset B_{V'}(x,c_0 r/\gamma_1^{p+1})$.

We will now use Theorem 3.1 with $x$ and $z$ as above, $t = \tau_{m-p}(z) > 0$ and
$$0 < \delta_1 = \frac{c_0 r}{L \gamma_1^{p+1}} < \delta_2 = C_1 L \rho_1^{p+1} < \omega(x)\;.$$
By Theorem 3.1, there exists a constant $K = K(\delta_1,\delta_2) > 0$ (depending on $\delta_1$ 
and $\delta_2$ which are  constants in our case) such that
$$\diam(B^u_t(z,\delta_1)) \geq \frac{1}{K\, G(z_{i_{m+1}})}\, \diam(B^u_t(z,\delta_2))
\geq \frac{1}{K\, G_0}\, \diam(B^u_t(z,\delta_2))\;.$$
However, using the above information about $\tcc$ and $\tcc'$, as in the proof of 
Proposition 3.3 in \cite{kn:St1}, 
one easily observes that $\cc' \supset B^u_t(z,\delta_1)$ and $\cc \subset B^u_t(z, \delta_2)$. Thus,
$\diam (\cc') \geq \frac{1}{K\, G_0}\, \diam(\cc)$. Combining the latter with (4.3) gives
$\diam (C[\ii']) \geq \frac{1}{L^2 K\, G_0}\, \diam(C[\ii])$.

This proves part (a) for $m > p$. Since there are only finitely many cylinders of length $\leq p$,
it follows immediately that there exists $\rho \in (0,1/(L^2KG_0)]$ which satisfies the requirements 
of part (a).

\ms

(b) This follows easily combining a simple modification of  the proof of 
Proposition 3.3(b) in \cite{kn:St1} with an argument similar to the above. 
We omit the details.
\endofproof

\bs

\noindent
{\it Proof of Theorem} 1.2. This now follows from the main result (Theorem 1.1) in \cite{kn:St1},
or rather from the proof of this theorem in section 5 in \cite{kn:St1}. What the latter assumes is
a local non-integrability condition (LNIC), uniformly Lispchitz local stable holonomy maps and 
the so called (see section 1 in \cite{kn:St1}) regular distortion along unstable manifolds. In our case
the flow is contact, so the condition (LNIC) follows from Proposition 6.2 in \cite{kn:St1}.
What concerns regular distortion along unstable manifolds, one should note that section 5 
in \cite{kn:St1} is only using a consequence of this property, namely the properties of cylinders 
described in Proposition 3.3 in \cite{kn:St1}. These properties are exactly the properties (a) and 
(b) in Theorem 4.2 above. Thus, under the assumptions of Theorem 1.2 above the argument 
from section 5 in \cite{kn:St1} applies and proves that the Ruelle transfer operators related to  
$f$  are eventually contracting for any Lipschitz real-valued function  $f$  on $U$.
\endofproof

\bs

\bs

\noindent
{\it University of Western Australia, Perth WA 6009, Australia\\
E-mail address: stoyanov@maths.uwa.edu.au}

\end{document}